\documentclass[12pt]{article}
\usepackage{amssymb}

\newenvironment{dfn}{\bigskip \noindent \bf Definition \rm}{\bigskip}
\newenvironment{rem}{\bigskip \noindent \bf Remark \rm}{\bigskip}
\newenvironment{example}{\bigskip \noindent \bf Example \rm}{\bigskip}
\newenvironment{proof}{\bigskip \noindent \bf Proof: \rm}{\bigskip}

\newenvironment{enumthm}{
\begin{enumerate}}{\end{enumerate}}
\newtheorem{thm}{Theorem}
\newtheorem{cor}[thm]{Corollary}
\newtheorem{lemma}[thm]{Lemma}
\newcommand{\qed}{\begin{flushright} \vspace{-1pc} $\square$
\end{flushright}}
\newfont{\knot}{orient scaled 1000}
\begin{document}
\begin{center}
{\Large\bf An Orientation-Sensitive Vassiliev Invariant \\[1ex]
for Virtual Knots} \\[3ex]
\large J\"org Sawollek\footnote{Fachbereich Mathematik, Universit\"at
Dortmund, 44221 Dortmund, Germany \\ {\em E-mail:\/}
sawollek@math.uni-dortmund.de \\ {\em WWW:\/}
http://www.mathematik.uni-dortmund.de/lsv/sawollek} \\[1ex]
March 13, 2002 (revised: August 30, 2002)
\end{center}
\vspace{4ex}

\begin{abstract}
It is an open question whether there are Vassiliev invariants that can
distinguish an oriented knot from its {\em inverse}, i.e., the knot with
the opposite orientation. In this article, an example is given for a
first order Vassiliev invariant that takes different values on a virtual
knot and its inverse. The Vassiliev invariant is derived from the Conway
polynomial for virtual knots. Furthermore, it is shown that the zeroth
order Vassiliev invariant coming from the Conway polynomial cannot
distinguish a virtual link from its inverse and that it vanishes for
virtual knots. \\[1ex]
\\[1ex]
\end{abstract}

\section*{Introduction}
In \cite{kauf} Kauffman defines virtual knot diagrams which are, in some
sense, a natural extension of classical knot diagrams. Surprisingly,
there can be found quite easily examples of virtual knots with
properties that are unknown for classical knots, e.g., a knot with
trivial Jones polynomial, see \cite{kauf}.

In this article, Vassiliev invariants for virtual knots and links are
investigated. They arise from the Conway polynomial for virtual links
defined in \cite{sawo}. It is shown that the lowest order coefficient of
the Conway polynomial vanishes on virtual knots and, furthermore, that
it is a zeroth order Vassiliev invariant for virtual links which cannot
detect a change of orientation. In contrast, two examples are given that
the corresponding Vassiliev invariant of order one distinguishes a
virtual knot from its inverse. For Vassiliev invariants of classical
knots it is still unknown whether they are orientation-sensitive or not.

\section{Virtual Links and Vassiliev Invariants}
\label{virtknots}

In the following, standard terminology from classical knot theory will
be used and some definitions that can be extended to virtual links in
the obvious way will not be given explicitly.

\begin{dfn}
A {\em virtual link diagram\/} is an oriented 4-valent
planar graph embedded in the plane with appropriate orientations of
edges and additional crossing information at each vertex (see
Fig. \ref{crossings}).
\begin{figure}[htb]
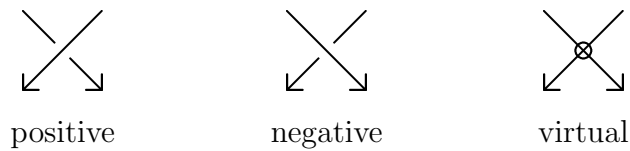

\vspace{-2ex}

\begin{center}
\begin{tabular}[t]{c@{\hspace{5em}}c@{\hspace{5em}}c}
\\[1ex]
{\knot p} & {\knot q} & {\knot r} \\[1ex]
positive & negative & virtual
\end{tabular}
\end{center}
\vspace{-1ex}

\caption{Crossing types}
\label{crossings}
\end{figure}
Denote the set of virtual link diagrams by $\cal
VD$. Two diagrams $D, D' \in {\cal VD}$ are called {\em equivalent\/} if
one can be transformed into the other by a finite sequence of extended
Reidemeister moves (see Fig. \ref{reidemeister})
\begin{figure}[htb]
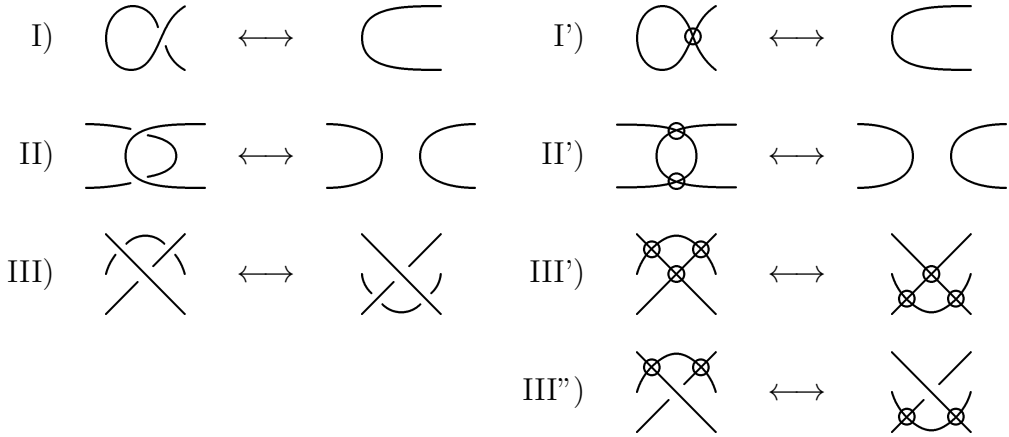

\begin{tabular}[t]{rccc}
I) & {\knot a} & $\longleftrightarrow$ & {\knot b} \\[2ex]
II) & {\knot d} & $\longleftrightarrow$ & {\knot e \hspace{0.5cm} f}
\\[2ex]
III) & {\knot g} & $\longleftrightarrow$ & {\knot h}
\end{tabular} \hspace{0em}
\begin{tabular}[t]{rccc}
I') & {\knot n} & $\longleftrightarrow$ & {\knot b} \\[2ex]
II') & {\knot m} & $\longleftrightarrow$ & {\knot e \hspace{0.5cm} f}
\\[2ex]
III') & {\knot k} & $\longleftrightarrow$ & {\knot l} \\[2ex]
III'') & {\knot i} & $\longleftrightarrow$ & {\knot j}
\end{tabular}
\vspace{2ex}

\caption{Extended Reidemeister moves}
\label{reidemeister}
\end{figure}
combined with orientation preserving homeomorphisms of the plane to
itself. A {\em virtual link\/} is an equivalence class of virtual link
diagrams. A virtual link with one component is called {\em virtual
knot}.
\end{dfn}

Though the notion of {\em components\/} of a virtual link, as in the
previous definition, comes from the idea to image the virtual link in
3-space as if the virtual crossings were classical ones, a virtual link
diagram does not correspond to an object in 3-space in the same way as
in classical knot theory, see \cite{sawo} for more details. Nevertheless
there exist several geometric interpretations for virtual links, see
\cite{cakasa}, \cite{kaka}, \cite{kauf}, \cite{sato}.

In the main part of this article, Vassiliev invariants for virtual
links, as introduced in \cite{kauf} (and not as in \cite{gpv}), will be
investigated. For an introduction to classical Vassiliev theory, see
\cite{BN} or \cite{CDL}.

\begin{dfn} A {\em singular virtual link diagram\/} is a virtual link
diagram that may contain vertices of an additional type called {\em
double points}, see Fig. \ref{double}.
\begin{figure}[htb]
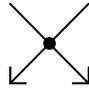

\begin{center}
\knot o
\end{center}
\vspace{-2ex}

\caption{A double point}
\label{double}
\end{figure} 
An equivalence relation is defined by adding the {\em rigid vertex
moves\/} depicted in Fig. \ref{rigid}
\begin{figure}[htb]
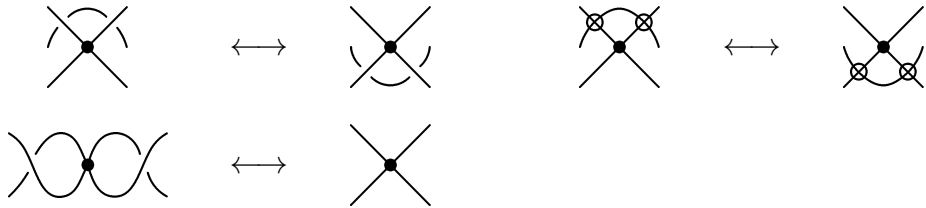

\begin{center}
\vspace{-2ex}

\begin{tabular}[t]{ccc@{\hspace{2cm}}ccc}
{\knot H} & $\quad \longleftrightarrow \quad$ & {\knot I} & {\knot F} &
$\quad \longleftrightarrow \quad$ & {\knot G} \\[2ex]
{\knot D} & $\quad \longleftrightarrow \quad$ & {\knot E}
\end{tabular}
\end{center}

\caption{Rigid vertex moves}
\label{rigid}
\end{figure} 
to the extended Reidemeister moves for virtual link diagrams.
\end{dfn}

\begin{figure}[htb]
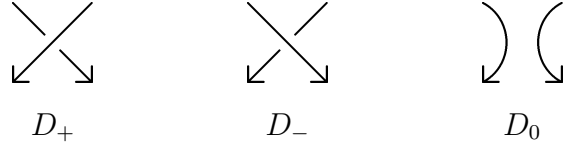

\begin{center}
\vspace{1ex}

\begin{tabular}[t]{c@{\hspace{5em}}c@{\hspace{5em}}c}
{\knot p} & {\knot q} & {\knot s} \\[1ex]
$D_+$ & $D_-$ & $D_0$
\end{tabular}
\end{center}
\vspace{-1ex}

\caption{Skein triple}
\label{skein}
\end{figure}

In the following, often diagrams are considered which are identical
except within a small disk where they differ as depicted in Fig.
\ref{skein}.
If $D$ is a singular virtual link diagram with a chosen double point
then replacing the double point with a positive or negative crossing
yields its positive and negative {\em resolution\/}, $D_+$ and $D_-$,
respectively. {\em Smoothing\/} the double point yields the
corresponding diagram $D_0$. Likewise, diagrams $D_{\varepsilon_1 \ldots
\varepsilon_k}$ with $\varepsilon_i \in \{ +, -, 0 \}$ can be defined
where $k$ chosen double points are replaced. The same notation is used
when crossings instead of double points are chosen and replaced.

\begin{dfn} A {\em Vassiliev invariant\/} $v$ is an invariant of
singular virtual link diagrams with values in an abelian group such that
\[ v(D_{\bullet}) \; = \; v(D_+) - v(D_-) \qquad \mbox{(Vassiliev
relation)} \]
where $D_{\bullet}$ denotes a diagram with a chosen double point and
$D_+$ and $D_-$ its positive and negative resolution, respectively. $v$
is said to be of {\em order\/} $\leq n$ if it vanishes on diagrams with
more than $n$ double points and it is said to be of {\em (exact)
order\/} $n$ if it is of order $\leq n$ but not of order $\leq n-1$.
\end{dfn}

\begin{rem}
Let $v$ denote an invariant of virtual links with values in an abelian
group. $v$ always can be extended to an invariant of singular virtual
links by demanding that the Vassiliev relation is fulfilled. $v$ is of
order zero iff $v(D_+) = v(D_-)$ holds for every virtual link diagram
$D$ and every crossing of $D$. $v$ is of order $\leq 1$ iff the equation
\[ v(D_{++}) - v(D_{+-}) - v(D_{-+}) + v(D_{--}) \; = \; 0 \]
holds for every virtual link diagram $D$ and every pair of crossings of
$D$.
\end{rem}

\section{Conway Polynomial and its Coefficients}
\label{conway}

The Conway polynomial for virtual links is derived from the normalized
$Z$-polynomial defined in \cite{sawo} which is an adaption of the Conway
polynomial for links in thickened surfaces that has been introduced by
Jaeger, Kauffman, and Saleur in \cite{jks}. In the following, the
construction of the $Z$-polynomial is recapitulated.

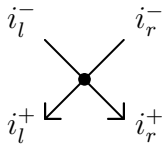
\begin{figure}[htb]
\begin{center}
\begin{picture}(1,1.5)
\put(0,0){\knot o}
\put(-0.5,1.2){$i_l^-$}
\put(-0.5,-0.2){$i_l^+$}
\put(1.2,1.2){$i_r^-$}
\put(1.2,-0.2){$i_r^+$}
\end{picture}
\end{center}

\caption{Half-edges at a vertex}
\label{vertex}
\end{figure}

Let $D$ be a virtual link diagram with $n \geq 1$ classical crossings
$c_1, \ldots, c_n$. Define
\[ M_+ := \left( \begin{array}{c@{\hspace{1em}}c}
1-x & -y \\
-xy^{-1} & 0
\end{array} \right) \, , \quad M_- := \left(
\begin{array}{c@{\hspace{1em}}c}
0 & -x^{-1}y \\
-y^{-1} & 1-x^{-1}
\end{array} \right) = M_+^{-1} \]
and $M := diag(M_+^{\varepsilon_1}, \ldots, M_+^{\varepsilon_n})$ where
$\varepsilon_i = \pm 1$ denotes the sign of $c_i$. Consider the graph
belonging to the virtual link diagram where the virtual crossings are
ignored, i.e., the graph consists of $n$ vertices $v_1, \ldots, v_n$
corresponding to the classical crossings and $2n$ edges corresponding to
the arcs connecting two classical crossings. Subdivide each edge into
two half-edges and label the four half-edges belonging to the vertex
$v_i$ as depicted in Fig. \ref{vertex}.

\noindent
The assignment
\[ (i,a) \mapsto (j,b) \quad \mbox{if $i_a^+$ and $j_b^-$ belong to the
same edge} \]
gives a permutation of the set $\{ 1, \ldots, n \} \times \{ l, r \}$.
Let $P$ denote the corresponding $2n \times 2n$ permutation matrix where
rows and columns are enumerated $(1,l)$, $(1,r)$, $(2,l)$, $(2,r)$,
\ldots, $(n,l)$, $(n,r)$, i.e., the $(i,a)$-th column of $P$ is the
$(j,b)$-th unit vector.

Define $Z_D(x,y) := (-1)^n \det(M-P)$ (observe that this yields the same
polynomial as the definition in \cite{sawo}). If $D$ has no classical
crossings then $Z_D(x,y) := 0$. See \cite{sawo} for example
calculations of $Z_D$.

\begin{thm}
$Z : {\cal VD} \to \mathbb Z[x^{\pm 1}, y^{\pm 1}]$ is an invariant of
virtual links up to multiplication by powers of $x^{\pm 1}$. To be more
precise, $Z$ is invariant with respect to all virtual Reidemeister moves
except moves of type I. Creating an additional (classical) crossing by
applying a Reidemeister move of type I yields a change of $Z$ by a
factor as shown in Fig. \ref{behaviour}.
\begin{figure}[htb]
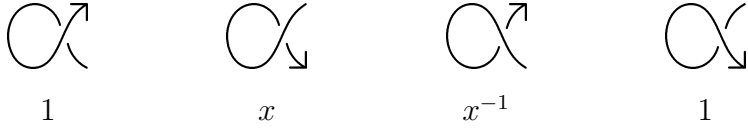

\begin{center}
\begin{tabular}[t]{c@{\hspace{4.5em}}c@{\hspace{4.5em}}c@{\hspace{4.5em}}c}
{\knot t} & {\knot u} & {\knot v} & {\knot w} \\[1ex]
$1$ & $x$ & $x^{-1}$ & $1$
\end{tabular}
\end{center}
\vspace{-2ex}

\caption{Behaviour of $Z$ with respect to Reidemeister moves of type I}
\label{behaviour}
\end{figure} 
\end{thm}

\begin{proof}
see \cite{sawo}
\qed
\end{proof}

\begin{rem}
Let $T := diag \left( \left( \begin{array}{cc}
0 & 1 \\
1 & 0
\end{array} \right), \ldots, \left( \begin{array}{cc}
0 & 1 \\
1 & 0
\end{array} \right) \right)$. The permutation corresponding to the
matrix $TP$ can be read off a virtual link diagram by walking along each
component of the link and writing down the numbers of the outgoing arcs.
This gives a cycle $\sigma_i$ for each link component and $\sigma_1
\ldots \sigma_r$ is the desired permutation. This describes an
alternative method to determine the matrix $P$.
\end{rem}

Define the {\em normalized polynomial\/} $\widetilde{Z}_D(x,y)$ as
follows. If $Z_D(x,y)$ is a non-vanishing polynomial and $N$ is the
lowest exponent in the variable $x$ then define
\[ \widetilde{Z}_D(x,y) \; := \; x^{-N} Z_D(x,y) \, . \]
Otherwise let $\widetilde{Z}_D(x,y) := Z_D(x,y) = 0$.

\begin{cor}
$\widetilde{Z} : {\cal VD} \to \mathbb Z[x, y^{\pm 1}]$ is an invariant
of virtual links.
\qed
\end{cor}

\begin{thm}
\label{props}
Let $D, D'$ denote virtual link diagrams and let $(D_+, D_-, D_0)$ be a
skein triple of virtual link diagrams. Then the following holds.
\begin{enumthm}
\item $Z_D(x,y) = 0$ if $D$ has no virtual crossings
\item $Z_{D \sqcup D'}(x,y) = Z_D(x,y) Z_{D'}(x,y)\;$ (disjoint
union)
\item $Z_D(1,y)$ does not depend on the over-under information of the
diagram's classical crossings
\item $Z_{D_+}(x,y) - x Z_{D_-}(x,y) = (1 - x) Z_{D_0}(x,y)$ (skein
relation of Conway type)
\end{enumthm}
\end{thm}

\begin{proof}
see \cite{sawo}
\qed
\end{proof}

\begin{rem}
Obviously, the parts a), b), c) of Theorem \ref{props} are valid for the
normalized polynomial $\widetilde{Z}$, too. Part d) is, in general, not
true for $\widetilde{Z}$ instead of $Z$ (simply consider a standard
diagram of a Hopf link with one classical crossing replaced by a virtual
crossing).
\end{rem}

Define the {\em Conway polynomial\/} $C_D(y,z) \in \mathbb Z[y^{\pm 1},
z]$ by expanding the normalized $Z$-polynomial in a Taylor series
\[ \widetilde{Z}_D(x,y) \; = \; \sum_k c_k (1-x)^k \]
and setting $z := 1-x$. Then $C_D(y,z)$ is an invariant of virtual links
which satisfies the Conway skein relation
\[ C_{D_+}(y,z) - C_{D_-}(y,z)  = z C_{D_0}(y,z) \]
"up to powers of $x = 1-z$", i.e., there exist integers $k_+, k_-, k_0$
such that
\[ (1-z)^{k_+} C_{D_+}(y,z) - (1-z)^{k_-} C_{D_-}(y,z)  = z (1-z)^{k_0}
C_{D_0}(y,z) \]
holds. Extend the Conway polynomial to singular virtual links via the
Vassiliev relation.

When the classical Conway polynomial $\nabla_D (z)$ is extended to
(classical) singular link diagrams by
\[ \nabla_{D_\bullet}(z) \; := \; \nabla_{D_+}(z) - \nabla_{D_-}(z) \]
then it can easily be seen that the coefficient $c_k$ of $\nabla_D(z) =
\sum c_k z^k$ is a Vassiliev invariant of order $\leq k$ since
\[ \nabla_{D_\bullet}(z) \; = \; z \nabla_{D_0}(z) \]
by the Conway skein relation. This does not work as easily for the
Conway polynomial $C_D(y,z)$, see \cite{talks}. For the purpose of this
article, it is enough to consider the Vassiliev invariants of orders
zero and one which are easy to handle.

\begin{lemma}
\label{formula}
Let $P$ denote the permutation matrix from the definition of $Z_D$, and
let $T = diag \left( \left( \begin{array}{cc}
0 & 1 \\
1 & 0
\end{array} \right), \ldots, \left( \begin{array}{cc}
0 & 1 \\
1 & 0
\end{array} \right) \right)$. Then
\[ c_0 \; = \; det(diag(y^{-1}, y, \ldots, y^{-1}, y) + TP). \]
\end{lemma}

\begin{proof}
Setting $x := 1$ in the definition of $Z_D$ gives:
\[ \begin{array}{lll}
c_0 & = & (-1)^n \cdot det \left( diag \left( \left( \begin{array}{cc}
0 & -y \\
-y^{-1} & 0
\end{array} \right) , \ldots , \left( \begin{array}{cc}
0 & -y \\
-y^{-1} & 0
\end{array} \right) \right) - P \right) \\[3ex]
& = & (-1)^n \cdot det(T) \cdot det( diag(-y^{-1}, -y, \ldots,
-y^{-1},-y) - TP)   \\[1ex]
& = & det( diag(y^{-1}, y, \ldots, y^{-1}, y) + TP )
\end{array} \]
\qed
\end{proof}

\begin{cor}
\label{or-inv}
$c_0$ is invariant with respect to a change of orientation.
\end{cor}

\begin{proof}
A change of orientation in a virtual link diagram corresponds to
replacing the permutation matrix $TP$ by the transposed matrix $(TP)^T$
which does not change the value of $c_0$ by Lemma \ref{formula}.
\qed
\end{proof}

In \cite{siwi} D. Silver and S. Williams show that their polynomial
invariant $\Delta_0$ is related to the $Z$-polynomial via an appropriate
change of variables (Proposition 3.1) and, furthermore, they investigate
the behaviour of $\Delta_0$ with respect to a change of orientation
(Corollary 5.2). Combining these results with Corollary \ref{or-inv}
immediately yields the following relation.

\begin{cor}
\label{symmetry}
For a virtual link with $d$ components,
\[ c_0(y) \; = \; (-1)^d c_0(y^{-1}) \]
where $c_0(y^{-1})$ denotes $c_0(y) = c_0$ with the variable $y$
replaced by $y^{-1}$.
\qed
\end{cor}

\begin{thm}
\label{vanish}
$c_0 = 0$ for virtual knots.
\end{thm}

\begin{proof}
The permutation corresponding to the matrix $TP$ can be read off a
virtual link diagram by walking along each component of the link and
writing down the numbers of the outgoing arcs. This gives a cycle for
each component. Therefore, if $D$ represents a virtual knot, then $TP$
is the permutation matrix of a cycle $\sigma$ and Lemma \ref{formula}
gives:
\[ \begin{array}{lll}
c_0 & = & det( diag(y^{-1}, y, \ldots, y^{-1}, y) + TP ) \\[1ex]
& = & sign(id) \cdot (y^{-1} \cdot y) \cdot \ldots \cdot (y^{-1} \cdot
y) + sign( \sigma ) \cdot 1 \cdot \ldots \cdot 1 \\[1ex]
& = & 1 + (-1)^{2n-1} \\[1ex]
& = & 0
\end{array} \]
\qed
\end{proof}

\begin{thm}
\label{vasinv}
$c_0$ is a Vassiliev invariant for virtual links, and $c_1$ is a
Vassiliev invariant for virtual knots.
\end{thm}

\begin{proof}
In the following, let $c_i^{\varepsilon}$ with $\varepsilon \in \{ +, -,
0 \}$ be an abbreviation for $c_i(D_{\varepsilon})$ and likewise
$c_i^{\varepsilon_1 \varepsilon_2} = c_i(D_{\varepsilon_1
\varepsilon_2})$ for a given virtual link diagram $D$ and crossings
chosen.

Since $c_0 = Z_D(1,y)$ it is clear from Theorem \ref{props} that $c_0$
is a Vassiliev invariant of order zero, i.e., $c_0^+ = c_0^-$.

For a virtual knot, $c_0$ vanishes by Theorem \ref{vanish} and therefore
a comparison of coefficients in the Conway skein relation
\[ (c_0^+ +(c_1^+ -k^+c_0^+)z + z^2(\ldots)) - (c_0^- +(c_1^-
-k^-c_0^-)z + z^2(\ldots)) = c_0^0z + z^2(\ldots) \]
immediately yields $c_1^+ - c_1^- = c_0^0$. Thus
\[ c_1^{++} - c_1^{+-} - c_1^{-+} + c_1^{--}  \; = \; c_0^{+0} -
c_0^{-0} \; = \; 0 \]
because $c_0$ is of order zero. This shows that $c_1$ is a Vassiliev
invariant of order one.
\qed
\end{proof}

\begin{rem}
Though formulated only for virtual links, Corollaries \ref{or-inv},
\ref{symmetry} and Theorems \ref{vanish}, \ref{vasinv} are valid for
singular virtual links, too, because of the Vassiliev relation.
\end{rem}

\section{Examples}
\label{examples}

\begin{example}
In general, $c_1$ is not a Vassiliev invariant of order one for virtual
links with more than one component. A counterexample is depicted in Fig.
\ref{notvas}.
\begin{figure}[htb]
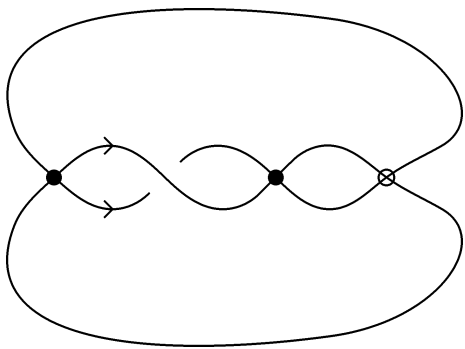

\begin{center}
\knot y
\end{center}

\caption{Singular virtual link with non-vanishing $c_1$}
\label{notvas}
\end{figure} 
$c_1$ does not vanish for the singular virtual link diagram $D$ shown in
Fig. \ref{notvas} though $D$ has two double points:
\[ c_1(D) \; = \; y+2+y^{-1} \]
Indeed, it is easy to construct examples which show that $c_1$ is
not a Vassiliev invariant of any order.
\end{example}

\begin{example}
In contrast to $c_0$, $c_1$ is orientation-sensitive. An example is
depicted in Fig. \ref{noninv}.
\begin{figure}[htb]
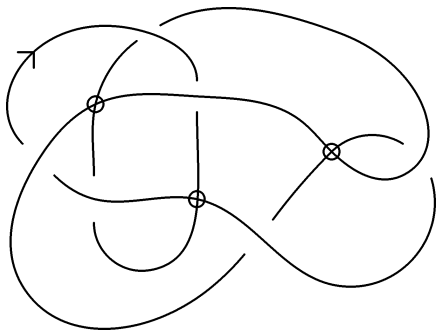

\begin{center}
\knot x
\end{center}

\caption{Non-invertible virtual knot}
\label{noninv}
\end{figure} 
Let $D$ be the diagram with the orientation indicated in Fig.
\ref{noninv} and let $D^*$ denote the diagram with the opposite
orientation. Then
\[ c_1(D) \; = \; -y^2-y+1+y^{-1} \qquad \mbox{and} \qquad c_1(D^*) \; =
\; y+1-y^{-1}-y^{-2}. \]
Chirality is detected, too. Let $\overline{D}$ denote the {\em mirror
diagram\/} of $D$, i.e., every classical crossing is changed from
positive to negative and vice versa. Then
\[ c_1(\overline{D}) \; = \; -y-1+y^{-1}+y^{-2} \qquad \mbox{and} \qquad
c_1(\overline{D}^*) \; = \; y^2+y-1-y^{-1}. \]
The Jones polynomial of $D$ (see \cite{kauf}) is non-trivial and it
takes different values on $D$ and $\overline{D}$. But from the
definition via the bracket polynomial it is clear that the Jones
polynomial cannot distinguish a virtual link from its inverse.
\end{example}

\begin{example}
Let $D$ be the diagram with the orientation indicated in Fig.
\ref{kauff}
\begin{figure}[htb]
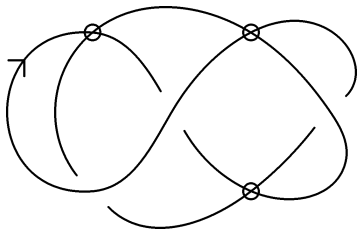

\begin{center}
\knot J
\end{center}

\caption{Kauffman's example}
\label{kauff}
\end{figure} 
and let $D^*$ denote the diagram with the opposite orientation. $D$ is
Kauffman's example of a virtual knot with trivial Jones polynomial and
trivial knot group (see \cite{kauf}). Again, $c_1$ distinguishes $D$ and
$D^*$. It happens that the values are the same as in the previous
example:
\[ c_1(D) \; = \; -y^2-y+1+y^{-1} \qquad \mbox{and} \qquad c_1(D^*) \; =
\; y+1-y^{-1}-y^{-2}. \]
\end{example}

\begin{rem}
In the same way as for classical links, quantum invariants can be
defined for virtual links, see \cite{gra}. It is well-known (and can be
shown analogously for virtual links) that quantum invariants cannot
detect a change of orientation (see \cite{kas} or \cite{tur}, for
example) and therefore $c_1$ is a Vassiliev invariant of virtual knots
which is not a function of quantum invariants. By a result of Gra\~{n}a
(\cite{gra}), this implies that $c_1$ neither is a function of quandle
cocycle invariants as defined in \cite{cjkls}.
\end{rem}

\section*{Acknowledgement}
The author would like to thank Dan Silver for valuable discussions.


\begin{thebibliography}{123}

\bibitem{BN} D. Bar-Natan, On the Vassiliev knot invariants, {\em
Topology\/} {\bf 34} (1995), 423--472.

\bibitem{cakasa} S. Carter, S. Kamada, and M. Saito, Stable equivalence
of knots on surfaces and virtual knot cobordisms, Preprint 2000. \\
{\em http://front.math.ucdavis.edu/math.GT/0008118}

\bibitem{cjkls} S. Carter, D. Jelsovsky, S. Kamada, L. Langford, M.
Saito, Quandle Cohomology and State-sum Invariants of Knotted Curves and
Surfaces, Preprint 1999. {\em
http://front.math.ucdavis.edu/math.GT/9903135}

\bibitem{CDL} S. Chmutov, S. Duzhin, and S. Lando, Vassiliev knot
invariants I. Introduction, {\em Adv. Soviet Math.\/} {\bf 21} (1994),
117--126.

\bibitem{gpv} M. Goussarov, M. Polyak, and O. Viro, Finite type
invariants of classical and virtual knots, {\em Topology\/} {\bf 39}
(2000), 1045--1068.

\bibitem{gra} M. Gra\~{n}a, Quandle knot invariants are quantum knot
invariants, {\em J. Knot Theory Ramifications\/} {\bf 11} (2002),
673--681.

\bibitem{jks} F. Jaeger, L.H. Kauffman, and H. Saleur, The Conway
Polynomial in $R^3$ and in Thickened Surfaces: A New Determinant
Formulation, {\em J. Combin. Theory Ser. B\/} {\bf 61} (1994), 237--259.

\bibitem{kaka} N. Kamada, and S. Kamada, Abstract link diagrams and
virtual knots, {\em J. Knot Theory Ramifications\/} {\bf 9} (2000),
93--106.

\bibitem{kas} C. Kassel, {\em Quantum Groups}, Springer, New York
(1995).

\bibitem{kauf} L.H. Kauffman, Virtual Knot Theory, {\em European J.
Combin.\/} {\bf 20} (1999), 663--690.

\bibitem{sato} S. Satoh, Virtual knot presentation of ribbon
torus-knots, {\em J. Knot Theory Ramifications\/} {\bf 9} (2000),
531--542.

\bibitem{sawo} J. Sawollek, On Alexander-Conway Polynomials for Virtual
Knots and Links, Preprint 1999. {\em
http://front.math.ucdavis.edu/math.GT/9912173}

\bibitem{talks} J. Sawollek, Talks given at Siegen in January 2001 and
at Bochum in April 2001, see \\
{\em http://www.mathematik.uni-dortmund.de/lsv/sawollek/bodoto01.html}

\bibitem{siwi} D. Silver, and S. Williams, Polynomial invariants of
virtual links, Preprint 2002. {\em http:// \\
www.southalabama.edu/mathstat/personal\_pages/williams/poly.pdf}

\bibitem{tur} V. Turaev, {\em Quantum Invariants of Knots and
3-Manifolds}, de Gruyter, Berlin (1994).

\end{thebibliography}
\end{document}